\documentclass[11pt]{article}
\usepackage{amsfonts}
\usepackage{color}

\usepackage{amsmath, amsthm, amsfonts, amssymb}

\newtheorem{lemma}{Lemma}
\newtheorem{theorem}{Theorem}
\newtheorem{proposition}{Proposition}
\newtheorem{remark}{Remark}
\newtheorem{corollary}{Corollary}

\def\btheorem{\begin{theorem}\sl{}\def\etheorem{\end{theorem}}}
\def\be{\begin{equation}}\def\ee{\end{equation}}
\def\bee{\begin{equation*}}\def\eee{\end{equation*}}
\def\beqlb{\begin{eqnarray}}\def\eeqlb{\end{eqnarray}}
\def\beqnn{\begin{eqnarray*}}\def\eeqnn{\end{eqnarray*}}
\def\bp{\begin{proposition}}\def\ep{\end{proposition}}

\def\qed{\hfill$\Box$\medskip}
\def\<{\langle}\def\>{\rangle}

\begin{document}
\centerline{\LARGE\bf Strong solutions of a class of SDEs with jumps
\footnote{Supported by NSFC (No. 10721091)}}

\bigskip

\centerline{\large \bf Juan Zhao\footnote{ \textit{E-mail address:}
{zhaojuan@mail.bnu.edu.cn }}}

\medskip

\centerline{{\it \footnotesize School of Mathematical Sciences,
Beijing Normal University,}}

\centerline{{\it \footnotesize Beijing 100875, People's Republic of
China}}

\medskip

\centerline{{\it \footnotesize Finalized October, 2008}}

\bigskip\bigskip

\noindent\rule[0.1cm]{16cm}{0.01cm}\\
\noindent{\bf Abstract}\\
\begin{small}

We study a class of stochastic integral equations with jumps under
non-Lipschitz conditions. We use the method of Euler approximations
to obtain the existence of the solution and give some sufficient
conditions for the strong uniqueness.

\end{small}

\bigskip
\begin{footnotesize}


\\\noindent\rule{16cm}{0.01cm}
\end{footnotesize}

\bigskip
\noindent{\bf 1. Introduction}\\
\setcounter{section}{1}
\setcounter{equation}{0}\setcounter{proposition}{0}\setcounter{theorem}{0}\setcounter{remark}{0}
\numberwithin{equation}{section} \numberwithin{proposition}{section}
\numberwithin{theorem}{section} \numberwithin{remark}{section}
\numberwithin{corollary}{section}

Modeling interest rate fluctuations is one of the major concerns of
both practitioners and academics. There are many prominent interest
rate models such as Vasicek model and Cox-Ingersoll-Ross model, see
Lamberton and Lapeyre (1996) for more details. Suppose that
$\{B(t)\}$ is a Brownian motion and $\{b(t)\}$ is a non-negative
measurable stochastic process. Let $\beta<0$ be a constant and
$\sigma$ be a $1/2$-H\"{o}lder continuous function on $\mathbb{R}_+$
vanishing at the origin. Deelstra and Delbaen (1995) introduced the
so-called extended CIR model $x(t)$ which is the solution of the
stochastic differential equation
 \beqnn dx(t)=(b(t)+\beta
x(t))dt+\sigma(x(t))dB(t)
 \eeqnn
with $x(0)\geq 0$. Deelstra and Delbaen (1998) used the method of
Euler approximations to prove the existence of the above stochastic
equation. In this paper, we extend the model by considering some
stochastic equations with jumps.

We consider a class of stochastic processes for the purpose of
modeling interest rates. Suppose that $U$ is a separable and
complete metric space. Let $\mu(du)$ be a $\sigma$-finite measure on
$U$. Let $(\Omega, \mathcal {F}, (\mathcal{F}_{t})_{t\geq0}, P)$ be
a filtered probability space satisfying the usual hypotheses. Let
$\{B(t)\}$ be a $(\mathcal {F}_{t})$-Brownian motion and let
$\{p(t)\}$ be a $(\mathcal {F}_{t})$-Poisson point process on $U$
with characteristic measure $\mu(du)$. Suppose that $\{B(t)\}$ and
$\{p(t)\}$ are independent of each other. Let $\{b(t)\}$ be a
non-negative measurable and adapted process and let $\{N(ds, du)\}$
be the Poisson random measure associated with $\{p(t)\}$. Suppose
that
\begin{itemize}
\item[(i)]$\beta<0$ is a constant and $x\mapsto \sigma(x)$
is a continuous function on $\mathbb{R}$ satisfying $\sigma(x)=0$
for $x\leq 0$;

\item[(ii)]$(x, u)\mapsto g(x, u)$ is a Borel function on
$\mathbb{R}\times U$ such that $g(x, u)+x\geq0$ for $x>0$ and $g(x,
u)=0$ for $x\leq 0$. \end{itemize} Given a non-negative
$\mathcal{F}_{0}$ measurable random variable $x(0)$, we consider the
following stochastic integral equation \beqlb\label{1.1}
x(t)=x(0)+\int_{0}^{t}(b(s)+\beta
x(s))ds+\int_{0}^{t}\sigma(x(s))dB(s)+\int_{0}^{t}\int_{U}g(x(s-),u)\tilde{N}(ds,
du) \eeqlb with $\tilde{N}(ds, du)=N(ds, du)-ds\mu(du).$ We are
interested in the existence and uniqueness of the solution for the
above stochastic equation. The coefficients of (\ref{1.1}) we are
considering are non-Lipschitz. Many authors studied the stochastic
equations which are closely related to the above equation. In
particular, Dawson and Li (2006, pp.1122-1131) gave a
characterization of continuous state branching processes with
immigration as strong solutions of some stochastic integral
equations. They used the tightness and the Skorokhod representation
to obtain the existence. Fu and Li (2008) studied a more general
class of stochastic equations with jumps. Under very weak
conditions, they established the existence and uniqueness of strong
solutions of those equations. The present work differs from that of
Fu and Li (2008) in that our drift term is given by a stochastic
process.

The remainder of the paper is organized as follows. In next section,
we state some results on the pathwise uniqueness of solutions to
(\ref{1.1}). In section 3, we discuss the Euler scheme for the
equation and show that the approximating solution converges in
$L^{1}$-supnorm towards the solution of (\ref{1.1}). Some criteria
on the existence and uniqueness of strong solutions are established
in the last section.

For some preliminary results concerning the stochastic differential
equations with jumps, the reader is referred to Bass (2004). We
refer to Ikeda and Watanabe (1989) and Protter (2004) for the theory
of stochastic analysis.

\bigskip

\noindent{\bf 2. Pathwise uniqueness}\\
\setcounter{section}{2}
\setcounter{equation}{0}\setcounter{proposition}{0}\setcounter{theorem}{0}\setcounter{remark}{0}
\numberwithin{equation}{section} \numberwithin{proposition}{section}
\numberwithin{theorem}{section} \numberwithin{remark}{section}
\numberwithin{corollary}{section}

In this section, we give some results on stochastic equations and on
the pathwise uniqueness of solutions to (\ref{1.1}). Because these
results can be obtained using essentially the same arguments as the
corresponding results of Fu and Li (2008), we omitted their proofs
here. Since the coefficients of (\ref{1.1}) satisfy the above
conditions, we have the following proposition.

\bp\label{p1} If $\{x(t)\}$ satisfies (\ref{1.1}) and $P\{x(0)\geq
0\}=1$, then $P\{x(t)\geq0 ~~for~ all ~~t\geq0\}=1$. \ep

In the sequel, we shall always assume the initial variable $x(0)$ is
non-negative, so Proposition \ref{p1} implies that any solution of
(\ref{1.1}) is non-negative. Then we can assume the ingredients are
defined only for $x\geq0$. In addition, for the convenience of the
statements of the results, we introduce the following conditions.

\begin{enumerate}

\item [(2.a)]The measurable and adapted process $b(\cdot)$ satisfying $\int_{0}^{t}Eb(s)ds<\infty$ for all
$t\geq0$;

\item [(2.b)]There is a constant $K\geq 0$ such that $\sigma^{2}(x)+\int_{U}\sup_{0\leq y\leq x}g^{2}(x,u)\mu(du)\leq K(1+x)$ for all $x\geq 0$;

\item [(2.c)]For every fixed $u\in U$, the function $x\mapsto g(x,u)$ is
non-decreasing, and for each integer $m\geq1$, there is a
non-negative and non-decreasing function $z\mapsto\rho_{m}(z)$ on
$\mathbb{R}_{+}$ so that $\int_{0+}\rho_{m}^{-2}(z)dz=\infty$ and
$|\sigma(x)-\sigma(y)|^{2}+\int_{U}[|l(x,y;u)|\wedge
l^{2}(x,y;u)]\mu(du)\leq \rho_{m}^{2}(|x-y|)$ for all $0\leq x,
y\leq m$, where $l(x,y;u)=g(x,u)-g(y,u)$;

\item [(2.d)]For every fixed $u\in U$, the function $x\mapsto g(x,u)$ is non-decreasing,
and for each integer $m\geq1$, there is a non-negative and
non-decreasing function $z\mapsto\rho_{m}(z)$ on $\mathbb{R}_{+}$ so
that $\int_{0+}\rho_{m}^{-2}(z)dz=\infty,~ |\sigma(x)-\sigma(y)|\leq
\rho_{m}(|x-y|)$ and $|g(x,u)-g(y,u)|\leq \rho_{m}(|x-y|)f_{m}(u)$
for all $0\leq x,y\leq m$ and $u\in U$, where $u\mapsto f_{m}(u)$ is
a non-negative function on $U$ satisfying $\int_{U}[f_{m}(u)\wedge
f_{m}^{2}(u)]\mu(du)<\infty$.
\end{enumerate}

We close this section with two theorems on the pathwise uniqueness
of solutions to (\ref{1.1}).

\btheorem\label{t1}Suppose that conditions (2.a, b, c) hold. Then
the pathwise uniqueness of solution holds for (\ref{1.1}). \etheorem

\btheorem\label{t2}Suppose that conditions (2.a, b, d) hold. Then
the pathwise uniqueness for (\ref{1.1}) holds. \etheorem

\noindent{\bf 3. Existence}\\
\setcounter{section}{3}
\setcounter{equation}{0}\setcounter{proposition}{0}\setcounter{theorem}{0}\setcounter{remark}{0}
\numberwithin{equation}{section} \numberwithin{proposition}{section}
\numberwithin{theorem}{section} \numberwithin{remark}{section}
\numberwithin{corollary}{section}

In this section, we prove a strong convergence of the Euler
approximations of the equation (\ref{1.1}), giving a construction of
the solution. A similar analysis was carried out in Yamada (1976,
1978) for continuous type equations, in Fu (2007, pp. 30-36) and Fu
and Li (2008) for two classes of jump-type equations. We refer the
reader to Mao et al. (2006, 2007) for recent results on related
topics.

For a fixed time $T>0$, we remark that $\int_{0}^{T}Eb(s)ds<\infty$.
Let us define the function $\gamma:
\mathbb{R_{+}}\rightarrow\mathbb{R_{+}}$ by \beqnn
\gamma(\nu)=\sup_{0\leq s\leq t\leq s+\nu\leq
T}\int_{s}^{t}Eb(u)du,~~ \nu\geq0. \eeqnn Since the function
$t\mapsto Eb(t)$ is integrable over the interval $[0,T]$, we have
that $\gamma(\nu)$ converges to zero as $\nu$ tends to zero.

We divide the interval $[0,T]$, known as the Euler discretization
method. For each $n\geq 1$, we take a subdivision
$$0=t_{0}^{n}\leq t_{1}^{n}\leq\cdots\leq t_{N_{n}}^{n}=T$$ and
denote this net by $\Delta_{n}$. For notational use, we drop the
index $n$ of the discretization times and write $N$ instead of
$N_{n}$.

The mesh of the net is defined as $\|\Delta_{n}\|=\sup_{1\leq k\leq
N}|t_{k}-t_{k-1}|$. We are working with a sequence of nets
$(\Delta_{n})_{n}$ such that the meshes are tending to zero. There
is no need to suppose that $\Delta_{n}\subset\Delta_{n+1}$.

The solutions of (\ref{1.1}) turns out to be non-negative but the
approximations we will need may take negative values. We therefore
put $\sigma'(x)=\sigma(x)I_{\{x\geq 0\}}$ and
$g'(x,u)=g(x,u)I_{\{x\geq 0\}}$. Note that $\sigma'(\cdot)$ and
$g'(\cdot, \cdot)$ also satisfy conditions (2.b, c, d).

If we are working with the net $\Delta_{n}$, we look at
$x_{\Delta_{n}}(t)$, which we denote by $x_{n}(t)$. We put
$x_{n}(0)=x(0)$. Let $\eta_{n}(t)=\sum_{k=0}^{k=N-1}t_{k}I_{[t_{k},
t_{k+1})}(t)$, we define a process $\{x_{n}(t)\}$
by\beqlb\label{3.1} x_{n}(t)&=&x(0)+\int_{0}^{t}(b(s)+\beta
x_{n}(\eta_{n}(s)))ds+\int_{0}^{t}\sigma'(x_{n}(\eta_{n}(s)))dB(s)\nonumber\\&&+\int_{0}^{t}\int_{U}g'(x_{n}(\eta_{n}(s)-),u)\tilde{N}(ds,du).\eeqlb
This is called an Euler approximation of (\ref{1.1}).

In the next conclusions, we need the following conditions:
\begin{enumerate}

\item [(3.a)]For every fixed $u\in U$, the function $x\mapsto g(x,u)$ is
non-decreasing, and for each integer $m\geq1$, there is a
non-negative and non-decreasing function $z\mapsto\rho_{m}(z)$ on
$\mathbb{R}_{+}$ so that $\int_{0+}\rho_{m}^{-2}(z)dz=\infty,
~z\mapsto \rho_{m}^{2}(z)$ is concave and
$|\sigma(x)-\sigma(y)|^{2}+\int_{U}l^{2}(x,y;u)\mu(du)\leq
\rho_{m}^{2}(|x-y|)$ for all $0\leq x, y\leq m$, where
$l(x,y;u)=g(x,u)-g(y,u)$;

\item [(3.b)]For every fixed $u\in U$, the function $x\mapsto g(x,u)$ is non-decreasing,
and for each integer $m\geq1$, there is a non-negative and
non-decreasing function $z\mapsto\rho_{m}(z)$ on $\mathbb{R}_{+}$ so
that $\int_{0+}\rho_{m}^{-2}(z)dz=\infty,~z\mapsto \rho_{m}^{2}(z)$
is concave, $|\sigma(x)-\sigma(y)|\leq \rho_{m}(|x-y|)$ and
$|g(x,u)-g(y,u)|\leq \rho_{m}(|x-y|)f_{m}(u)$ for all $0\leq x,y\leq
m$ and $u\in U$, where $u\mapsto f_{m}(u)$ is a non-negative
function on $U$ satisfying $\int_{U}f_{m}^{2}(u)\mu(du)<\infty$.
\end{enumerate}
It is easy to show that $\sigma'(\cdot)$ and $g'(\cdot, \cdot)$ also
satisfy conditions (3.a, b).

\begin{remark}The functions $\rho(z)=\sqrt{z},~\rho(z)=z^{\frac{1}{2}}(\log\frac{1}{z})^{\frac{1}{2}},
~\rho(z)=z^{\frac{1}{2}}(\log\frac{1}{z})^{\frac{1}{2}}(\log\log\frac{1}{z})^{\frac{1}{2}},
\cdots$ satisfy conditions (3.a, b).\end{remark}

\btheorem\label{t3}Suppose that conditions (2.a, b) and (3.a) hold.
Then the Euler scheme (\ref{3.1}) with $t_{k}\leq t<t_{k+1}, ~k=0,
1, \cdots, N-1$ converges to the solution of (\ref{1.1}) in
$L^{1}$-supnorm. \etheorem

\begin{remark} If the intensity of the Poisson random measure is
zero and $\rho(z)=\sqrt{z}$, the results are degenerated to those of
Deelstra and Delbaen (1998).\end{remark}

Next, we prove $x_{n}(t)$ converges to the solution of (\ref{1.1})
in $L^{1}$-supnorm.

\bp\label{p4} Suppose that condition (2.a, b) hold. Then for all
$0\leq t\leq T$, there exist constants $G_{T}\geq 0$ and $H_{T}\geq
0$ such that the following hold: \beqlb\label{3.2}
E[|x_{n}(\eta_{n}(t))|]&\leq& G_{T};\\
\label{3.3}E[|x_{n}(t)|]&\leq& H_{T};\\
\label{3.4}E[|x_{n}(t)-x_{n}(\eta_{n}(t))|]&\leq&
\gamma(\|\Delta_{n}\|)-\beta
G_{T}\|\Delta_{n}\|+2\sqrt{K(G_{T}+1)\|\Delta_{n}\|}. \eeqlb

\ep

\noindent{\bf Proof.} From (\ref{3.1}), we obtain \beqnn
E[|x_{n}(\eta_{n}(t))|] &\leq&
E[x(0)]+\int_{0}^{t}Eb(s)ds+E^{\frac{1}{2}}[\int_{0}^{t}\sigma'(x_{n}(\eta_{n}(s)))dB(s)]^{2}\\
&&+|\beta|\int_{0}^{t}E[|x_{n}(\eta_{n}(s))|]ds+E^{\frac{1}{2}}[\int_{0}^{t}\int_{U}g'(x_{n}(\eta_{n}(s)-),u)\tilde{N}(ds,du)]^{2}\\
&\leq&
E[x(0)]+\int_{0}^{t}Eb(s)ds+|\beta|\int_{0}^{t}E[|x_{n}(\eta_{n}(s))|]ds\\
&& +2+E[\int_{0}^{t}\sigma'^{2}(x_{n}(\eta_{n}(s)))ds]+E[\int_{0}^{t}ds\int_{U}g'^{2}(x_{n}(\eta_{n}(s)-),u)\mu(du)]\\
&\leq&
E[x(0)]+\int_{0}^{t}Eb(s)ds+|\beta|\int_{0}^{t}E[|x_{n}(\eta_{n}(s))|]ds\\
&&
+2+Kt+K\int_{0}^{t}E[|x_{n}(\eta_{n}(s))|]ds\\
&\leq&
(E[x(0)]+\int_{0}^{t}Eb(s)ds+2+Kt)+(K-\beta)\int_{0}^{t}E[|x_{n}(\eta_{n}(s))|]ds.
\eeqnn The first and the third inequalities follow by Cauchy-Schwarz
inequality and condition (2.b) respectively. By Gronwall's lemma, we
get \beqnn E[|x_{n}(\eta_{n}(t))|]&\leq&
(E[x(0)]+\int_{0}^{t}Eb(s)ds+2+Kt)\exp\{(K-\beta)t\}\\&=:&G_{t}\leq
G_{T}. \eeqnn

After similar calculations, from (\ref{3.1}) and (\ref{3.2}), we get
\beqnn E[|x_{n}(t)|]&\leq&
E[x(0)]+\int_{0}^{t}Eb(s)ds+2+Kt+(K-\beta)tG_{t}\\&=:&H_{t}\leq
H_{T}.\eeqnn The above two bounds are independent of $n$ and $t$.

From (\ref{3.1}), (\ref{3.2}), (\ref{3.3}) and Cauchy-Schwarz
inequality, we get (\ref{3.4}) immediately. \qed

Given a function $f$ defined on a subset of $\mathbb{R}$, we note
\beqnn \Delta_{z}f(x)=f(x+z)-f(x)~ ~and~ ~
D_{z}f(x)=\Delta_{z}f(x)-f'(x)z\eeqnn if the right hand sides are
meaningful.

\bp\label{p5} Suppose that conditions (2.a, b) and (3.a) hold. Then
there exists a progressive process $\{y(t)\}$ such that the
following convergence hold: \beqlb\label{3.5}
\lim_{n\rightarrow\infty}\sup_{0\leq t\leq T}E[|x_{n}(t)-y(t)|]&=&0;\\
\label{3.6}\lim_{n\rightarrow\infty}\sup_{0\leq t\leq
T}E[|x_{n}(\eta_{n}(t))-y(t)|]&=&0. \eeqlb \ep

\noindent{\bf Proof.}~Let $\zeta(t)=x_{n}(t)-x_{n'}(t)$ for fixed
$n,~n'\geq 1$. Following from (\ref{3.1}), we get \beqlb\label{3.7}
\zeta(t)&=&\beta\int_{0}^{t}[x_{n}(\eta_{n}(s))-x_{n'}(\eta_{n'}(s))]ds
+\int_{0}^{t}[\sigma'(x_{n}(\eta_{n}(s)))-\sigma'(x_{n'}(\eta_{n'}(s)))]dB(s)\nonumber\\
&&+\int_{0}^{t}\int_{U}[g'(x_{n}(\eta_{n}(s)-),u)-g'(x_{n'}(\eta_{n'}(s)-),u)]\tilde{N}(ds,du).\eeqlb

Let $a_{0}=1$ and choose $a_{k}\to 0+$ decreasingly so that
$\int_{a_{k}}^{a_{k-1}}\rho_{m}^{-2}(z)dz=k$ for $k\geq 1$. Let
$z\mapsto\psi_{k}(z)$ be a non-negative continuous function on
$\mathbb{R}$ which has support in $(a_{k},a_{k-1})$ and satisfies
$\int_{a_{k}}^{a_{k-1}}\psi_{k}(z)dz=1$ and $0\leq
\psi_{k}(z)\leq2k^{-1}\rho_{m}^{-2}(z)$ for $a_{k}<z<a_{k-1}$. For
each $k\geq 1$ we define the non-negative and twice continuously
differentiable function \beqnn
\phi_{k}(z)=\int_{0}^{|z|}dy\int_{0}^{y}\psi_{k}(x)dx,~~
z\in\mathbb{R}.\eeqnn Clearly, the sequence $\{\phi_{k}\}$ satisfies
\begin{enumerate}

\item [(i)]$\phi_{k}(x)\rightarrow|x|$ non decreasingly as $k\rightarrow\infty$;

\item [(ii)]$0\leq \phi_{k}'(x)\leq 1$ for $x\geq 0$ and $-1\leq \phi_{k}'(x)\leq 0$ for $x\leq
0$.
\end{enumerate} Let
$\tau_{m}=\inf\{t\geq 0, |x_{n}(t)|\geq m~ or ~|x_{n'}(t)|\geq m\}$
for $m\geq1$. Applying It\^{o}'s formula, we get \beqlb\label{3.8}
\phi_{k}(\zeta(t\wedge \tau_{m})) &=&
\beta\int_{0}^{t\wedge \tau_{m}}\phi_{k}'(\zeta(s))[x_{n}(\eta_{n}(s))-x_{n'}(\eta_{n'}(s))]ds+mart. \nonumber\\
&&
+\frac{1}{2}\int_{0}^{t\wedge \tau_{m}}\phi_{k}''(\zeta(s))[\sigma'(x_{n}(\eta_{n}(s)))-\sigma'(x_{n'}(\eta_{n'}(s)))]^{2}ds\nonumber\\
&&
+\int_{0}^{t\wedge \tau_{m}}ds\int_{U}[D_{l(n,n';u)}\phi_{k}(\zeta(s-))]\mu(du)\nonumber\\
&=:& I_{1}(t\wedge \tau_{m})+mart.+I_{2}(t\wedge
\tau_{m})+I_{3}(t\wedge \tau_{m}), \eeqlb where
$l(n,n';u)=g'(x_{n}(\eta_{n}(s)-),u)-g'(x_{n'}(\eta_{n'}(s)-),u)$.

According to $\beta<0$ and $\{\phi_{k}\}$ satisfies property (ii),
we get
 \beqnn
&&I_{1}(t\wedge \tau_{m})\\
 &=&\beta\int_{0}^{t\wedge \tau_{m}}\phi_{k}'(\zeta(s))[x_{n}(\eta_{n}(s))-x_{n}(s)]ds+\beta\int_{0}^{t\wedge \tau_{m}}\phi_{k}'(\zeta(s))[x_{n}(s)-x_{n'}(s)]ds\\
&&+\beta\int_{0}^{t\wedge \tau_{m}}\phi_{k}'(\zeta(s))[x_{n'}(s)-x_{n'}(\eta_{n'}(s))]ds\\
&\leq& \beta\int_{0}^{t\wedge
\tau_{m}}\phi_{k}'(\zeta(s))[x_{n}(\eta_{n}(s))-x_{n}(s)]ds+\beta\int_{0}^{t\wedge
\tau_{m}}\phi_{k}'(\zeta(s))[x_{n'}(s)-x_{n'}(\eta_{n'}(s))]ds.\eeqnn
Consequently,
 \beqnn
&&E[I_{1}(t\wedge \tau_{m})]\\
 &\leq&
|\beta|E[\int_{0}^{t\wedge \tau_{m}}|x_{n}(\eta_{n}(s))-x_{n}(s)|ds]
+|\beta|E[\int_{0}^{t\wedge \tau_{m}}|x_{n'}(s)-x_{n'}(\eta_{n'}(s))|ds]\\
 &=:&
|\beta|A(n,m,t)+|\beta|A(n',m,t),
 \eeqnn
where $A(n,m,t)=E[\int_{0}^{t\wedge
\tau_{m}}|x_{n}(\eta_{n}(s))-x_{n}(s)|ds]$.

Since $\int_{a_{k}}^{a_{k-1}}\rho_{m}^{-2}(z)dz=k$ and by the
monotonicity of $z\mapsto \rho_{m}(z)$, we have
$k^{-1}\rho_{m}^{-2}(a_{k})\leq 2$. Note that $0\leq
\phi_{k}''(z)=\psi_{k}(|z|)\leq
2k^{-1}\rho_{m}^{-2}(|z|)I_{(a_{k},a_{k-1})}(|z|)\leq2k^{-1}\rho_{m}^{-2}(|a_{k}|)\leq
4$ and $\sigma'(x)$ satisfies condition (3.a), we have
 \beqnn
E[I_{2}(t\wedge \tau_{m})]&\leq&\frac{3}{2}E[\int_{0}^{t\wedge
\tau_{m}}\phi_{k}''(\zeta(s))(\sigma'(x_{n}(\eta_{n}(s)))-\sigma'(x_{n}(s)))^{2}ds]\\
&&+\frac{3}{2}E[\int_{0}^{t\wedge \tau_{m}}\phi_{k}''(\zeta(s))(\sigma'(x_{n}(s))-\sigma'(x_{n'}(s)))^{2}ds]\\
&&+\frac{3}{2}E[\int_{0}^{t\wedge \tau_{m}}\phi_{k}''(\zeta(s))(\sigma'(x_{n'}(s))-\sigma'(x_{n'}(\eta_{n'}(s))))^{2}ds]\\
&\leq&6E[\int_{0}^{t\wedge \tau_{m}}\rho_{m}^{2}(|x_{n}(\eta_{n}(s))-x_{n}(s)|)ds]+\frac{3t}{k}\\
&&+6E[\int_{0}^{t\wedge
\tau_{m}}\rho_{m}^{2}(|x_{n'}(s)-x_{n'}(\eta_{n'}(s))|)ds]\\
&=:&6B(n,m,t)+\frac{3t}{k}+6B(n',m,t), \eeqnn where
$B(n,m,t)=E[\int_{0}^{t\wedge
\tau_{m}}\rho_{m}^{2}(|x_{n}(\eta_{n}(s))-x_{n}(s)|)ds]$.

By Taylor's expansion and the definition of $\phi_{k}$, for all $h,
~\zeta\in \mathbb{R}$ it is easy to show that \beqlb\label{3.9}
D_{h}\phi_{k}(\zeta)&=&h^{2}\int_{0}^{1}\phi_{k}''(\zeta+th)(1-t)dt=h^{2}\int_{0}^{1}\psi_{k}(|\zeta+th|)(1-t)dt\nonumber\\
&\leq&
2k^{-1}h^{2}\int_{0}^{1}\rho_{m}^{-2}(|\zeta+th|)\{I_{(a_{k},a_{k-1})}(|\zeta+th|)\}(1-t)dt\nonumber\\&\leq&
2k^{-1}h^{2}\rho_{m}^{-2}(a_{k})\int_{0}^{1}(1-t)dt \leq
2h^{2}.\eeqlb Note also that $\zeta(s-)\neq\zeta(s)$ for at most
countably many $s\geq 0$. From (\ref{3.8}), (\ref{3.9}) and
$g'(x,u)$ satisfies condition (3.a), we have \beqnn E[I_{3}(t\wedge
\tau_{m})]&\leq&2E[\int_{0}^{t\wedge
\tau_{m}}ds\int_{U}(u_{1}+u_{2}+u_{3})^{2}\mu(du)]
\leq6\sum_{i=1}^{3}E[\int_{0}^{t\wedge
\tau_{m}}ds\int_{U}u_{i}^{2}\mu(du)]\\
&\leq&6E[\int_{0}^{t\wedge
\tau_{m}}\rho_{m}^{2}(|x_{n}(\eta_{n}(s))-x_{n}(s)|)ds]+6E[\int_{0}^{t\wedge
\tau_{m}}\rho_{m}^{2}(|x_{n}(s)-x_{n'}(s)|)ds]\\
&&+6E[\int_{0}^{t\wedge
\tau_{m}}\rho_{m}^{2}(|x_{n'}(s)-x_{n'}(\eta_{n'}(s))|)ds]\\
&=&6B(n,m,t)+6E[\int_{0}^{t\wedge
\tau_{m}}\rho_{m}^{2}(|\zeta(s)|)ds]+6B(n',m,t),\eeqnn where \beqnn
u_{1}&=&g'(x_{n}(\eta_{n}(s)-),u)-g'(x_{n}(s-),u),\\
u_{2}&=&g'(x_{n}(s-),u)-g'(x_{n'}(s-),u), \\
u_{3}&=&g'(x_{n'}(s-),u)-g'(x_{n'}(\eta_{n'}(s)-),u).\eeqnn

Consequently, \beqlb\label{3.10} E[\phi_{k}(\zeta(t\wedge
\tau_{m}))]&\leq& -\beta A(n,m,t)-\beta
A(n',m,t)+12B(n,m,t)+12B(n',m,t)\nonumber\\&&+\frac{3t}{k}+6E[\int_{0}^{t\wedge
\tau_{m}}\rho_{m}^{2}(|\zeta(s)|)ds].\eeqlb By Proposition \ref{p4},
the assumption on $z\mapsto \rho_{m}^{2}(z)$ and the dominated
convergence theorem, it is easy to see that \beqnn
\lim_{n\rightarrow\infty}A(n,m,t)=\lim_{n\rightarrow\infty}B(n,m,t)=0.\eeqnn
From the definition of $\phi_{k}(\cdot)$, we remark that $|z|\leq
a_{k-1}+\phi_{k}(z)$ for every $z\in \mathbb{R}$. For given $T\geq0$
and $\varepsilon>0$, we first take an integer $k_{0}\geq 1$ such
that $a_{k_{0}-1}+3T/k_{0}<\varepsilon/2$. Then we choose
sufficiently large $N=N(k_{0})\geq1$ so that $12B(n,m,t)-\beta
A(n,m,t)<\varepsilon/4$ for every $n\geq N$. By (\ref{3.10}), we
have \beqnn E[|\zeta(t\wedge \tau_{m})|]\leq
\varepsilon+6E[\int_{0}^{t\wedge
\tau_{m}}\rho_{m}^{2}(|\zeta(s)|)ds]\eeqnn for $0\leq t\leq T$ and
$n,n'\geq N$. Since $\zeta(s)<2m$ for all $0<s\leq \tau_{m}$, we
infer that $t\mapsto E[|\zeta(t\wedge \tau_{m})|]$ is locally
bounded. Then the concaveness of $z\mapsto \rho_{m}^{2}(z)$ implies
that \beqlb\label{3.11} E[|\zeta(t\wedge \tau_{m})|]&\leq&
\varepsilon+6E[\int_{0}^{t}\rho_{m}^{2}(|\zeta(s\wedge
\tau_{m})|)ds]\nonumber\\&\leq&
\varepsilon+6\int_{0}^{t}\rho_{m}^{2}(E(|\zeta(s\wedge
\tau_{m})|))ds\eeqlb for $0\leq t\leq T$ and $n,n'\geq N$. Let
\beqnn R_{n}(t)=\sup_{n'\geq n}\sup_{0\leq s\leq t}E[|x_{n}(s\wedge
\tau_{m})-x_{n'}(s\wedge \tau_{m})|], ~~n\geq 1, 0\leq t\leq
T.\eeqnn In view of (\ref{3.11}), the monotonicity of $z\mapsto
\rho_{m}^{2}(z)$ gives \beqnn R_{n}(t)\leq
\varepsilon+6\int_{0}^{t}\rho_{m}^{2}(R_{n}(s))ds, ~~n\geq N, ~0\leq
t\leq T.
 \eeqnn
 By letting $n \rightarrow \infty$ and $\varepsilon \rightarrow 0$ we obtain
  \beqnn
  \lim_{n\rightarrow \infty}R_{n}(t)\leq
  6\int_{0}^{t}\rho_{m}^{2}(\lim_{n\rightarrow\infty}R_{n}(s))ds,
  ~0\leq t\leq T.
  \eeqnn
  Thus $\lim_{n\rightarrow\infty}R_{n}(t)=0$ for every $0\leq t\leq T$. Since $\tau_{m}\to \infty$ as
$m\to\infty$ by Proposition \ref{p4}, now letting $m\to \infty$, it
is easy to find a progressive process $\{y(t)\}$ such that
(\ref{3.5}) holds. Moreover, by (\ref{3.4}) and (\ref{3.5}),
(\ref{3.6}) is also obtained.  \qed

\bp\label{p6} Suppose that conditions (2.a, b) and (3.a) hold. Then
there exists a c\`{a}dl\`{a}g process $\{x(t)\}$ such that: \beqnn
\lim_{n\rightarrow\infty}E[\sup_{0\leq t\leq
T}|x_{n}(t)-x(t)|]&=&0;\\
\lim_{n\rightarrow\infty}\sup_{0\leq t\leq
T}E[|x_{n}(\eta_{n}(t))-x(t)|]&=&0 \eeqnn hold. Moreover, $\{x(t)\}$
is a non-negative solution of (\ref{1.1}).\ep

\noindent{\bf Proof.}~Let $\tau_{m}=\inf\{t\geq 0, x_{n}(t)\geq m~
or ~x_{n'}(t)\geq m\}$ for $m\geq1$. Applying Doob's martingale
inequality to (\ref{3.7}), we get\beqnn E[\sup_{0\leq t\leq
T}|x_{n}(t\wedge\tau_{m})-x_{n'}(t\wedge\tau_{m})|]
&\leq&|\beta|\int_{0}^{T\wedge\tau_{m}}E[|x_{n}(\eta_{n}(s))-x_{n'}(\eta_{n'}(s))|]ds\\
&&+4E^{\frac{1}{2}}[\int_{0}^{T\wedge\tau_{m}}(\sigma'(x_{n}(\eta_{n}(s)))-\sigma'(x_{n'}(\eta_{n'}(s))))^{2}ds]\\
&&+4E^{\frac{1}{2}}[\int_{0}^{T\wedge\tau_{m}}ds\int_{U}l^{2}(n,n';u)\mu(du)].\eeqnn

Letting $m\to\infty$, by condition (3.a), Proposition \ref{p4},
\ref{p5} and dominated convergence theorem, we get\beqnn \lim_{n,
n'\rightarrow\infty}E[\sup_{0\leq t\leq T}|x_{n}(t)-x_{n'}(t)|]=0.
\eeqnn

Consequently, $\{y(t)\}$ has a c\`{a}dl\`{a}g modification
$\{x(t)\}$ satisfying the first equality. The second equality then
follows by Proposition \ref{p5}.

Next, we will show that \beqnn x(t)=x(0)+\int_{0}^{t}(b(s)+\beta
x(s))ds+\int_{0}^{t}\sigma'(x(s))dB(s)
+\int_{0}^{t}\int_{U}g'(x(s-),u)\tilde{N}(ds, du). \eeqnn Indeed,
from (\ref{3.1}) \beqnn &&E[\sup_{0\leq t\leq
T}|x(t)-x(0)-\int_{0}^{t}(b(s)+\beta
x(s))ds-\int_{0}^{t}\sigma'(x(s))dB(s)\\
&&-\int_{0}^{t}\int_{U}g'(x(s-),u)\tilde{N}(ds, du)|]\\
&=&E[\sup_{0\leq t\leq T}|x(t)-x_{n}(t)+\beta
\int_{0}^{t}(x_{n}(\eta_{n}(s))-x(s))ds+\int_{0}^{t}(\sigma'(x_{n}(\eta_{n}(s)))-\sigma'(x(s)))dB(s)\\
&&+\int_{0}^{t}\int_{U}(g'(x_{n}(\eta_{n}(s)-),u)-g'(x(s-),u))\tilde{N}(ds,
du)|]\eeqnn and the result follows by the triangular inequality,
Doob's martingale inequality and the previous calculations.

By Proposition \ref{p1} and the definitions of $\sigma'(x)$ and
$g'(x,u)$, we remark that $x(t)$ is a non-negative process.
Therefore, we can replace $\sigma'(x)$ and $g'(x,u)$ by $\sigma(x)$
and $g(x,u)$ respectively. Consequently, $x(t)$ satisfies \beqnn
x(t)=x(0)+\int_{0}^{t}(b(s)+\beta
x(s))ds+\int_{0}^{t}\sigma(x(s))dB(s)
+\int_{0}^{t}\int_{U}g(x(s-),u)\tilde{N}(ds, du). \eeqnn Then we
complete the proof.\qed

We prove that the Euler scheme (\ref{3.1}) converges to the unique
solution of (\ref{1.1}) in $L^{1}$-supnorm. The conclusion of
Theorem \ref{t3} holds immediately.

After similar analysis to the previous results, we have the
following theorem.

\btheorem\label{t4}Suppose that conditions (2.a, b) and (3.b) are
satisfied. Then the Euler scheme (\ref{3.1}) with $t_{k}\leq
t<t_{k+1}, ~k=0, 1, \cdots, N-1$ converges in $L^{1}$-supnorm
towards the solution of (\ref{1.1}). \etheorem

\bigskip
\noindent{\bf 4. Strong solutions}\\
\setcounter{section}{4}
\setcounter{equation}{0}\setcounter{proposition}{0}\setcounter{theorem}{0}\setcounter{remark}{0}
\numberwithin{equation}{section} \numberwithin{proposition}{section}
\numberwithin{theorem}{section} \numberwithin{remark}{section}
\numberwithin{corollary}{section}

In this section, we give some criteria on the existence and
uniqueness of the strong solution of equation (\ref{1.1}) and
illustrate a simple application of the results to stochastic
differential equations driven by one-sided L\'{e}vy processes.

\btheorem\label{t5}Suppose that conditions (2.a, b) and (3.a) are
satisfied. Then there exists a unique non-negative strong solution
to (\ref{1.1}). \etheorem

\noindent{\bf Proof}.~ By applying Theorem \ref{t3} we infer that
(\ref{1.1}) has a non-negative solution. In addition, the
ingredients of (\ref{1.1}) satisfy condition (2.c). Then the
pathwise uniqueness of the equation follows from Theorem \ref{t1}.
\qed

Based on the pathwise uniqueness stated in Theorem \ref{t2}, the
following result can be proved similarly as the above.

\btheorem\label{t6}Suppose that conditions (2.a, b) and (3.b) are
satisfied. Then there exists a unique non-negative strong solution
to (\ref{1.1}). \etheorem

At last, we give a simple application of Theorem \ref{t5} and
Theorem \ref{t6}. Now we consider stochastic equations driven by
one-sided L\'{e}vy processes. Let $\mu(dz)$ be a $\sigma$-finite
measure on $(0, \infty)$. We assume that
$\int_{0}^{\infty}z^{2}\mu(dz)<\infty.$ Let $\{B(t)\}$ be a standard
$(\mathcal {F}_{t})$-Brownian motion. Let $\{z(t)\}$ be a $(\mathcal
{F}_{t})$-L\'{e}vy process with exponent $
u\mapsto\int_{0}^{\infty}(e^{iuz}-1-iuz)\mu(dz)$. Therefore
$\{z(t)\}$ is centered. Suppose that those processes are independent
of each other. In addition, suppose that $\beta<0$ is a real
constant and
\begin{enumerate}
\item[(i)]a measurable and adapted process $b:\Omega \times \mathbb{R}_{+}\mapsto\mathbb{R}_{+}$
 satisfying $\int_{0}^{t}Eb(s)ds<\infty$ for all $t\geq0$.

\item[(ii)]  $x\mapsto \sigma(x)$ is a continuous
function on $\mathbb{R}_{+}$ satisfying $\sigma(0)=0$;

\item [(iii)]$x\mapsto \phi(x)$ is a continuous non-negative function on $\mathbb{R}_{+}$ satisfying $\phi(0)=0$.
\end{enumerate}
We assume the following condition on the
ingredients:\begin{enumerate}\item [(4.a)] The function
$x\mapsto\phi(x)$ is non-decreasing and for each $m\geq 1$ there is
a constant $K_{m}\geq 0$ so that \beqnn
|\sigma(x)-\sigma(y)|^{2}+|\phi(x)-\phi(y)|^{2}\leq K_{m}|x-y|
\eeqnn for all $0\leq x,y \leq m$.\end{enumerate}

\btheorem\label{t7} Under condition (4.a), there is a unique
non-negative strong solution to \beqlb\label{4.1}
dx(t)=\sigma(x(t))dB(t)+\phi(x(t-))dz(t)+(b(t)+\beta
x(t))dt.\eeqlb\etheorem

{\noindent\bf Proof.}~By the general result on L\'{e}vy-It\^{o}
decompositions, see, e.g., Sato (1999, p.120, Theorem 19.2), we have
\beqnn z(t)=\int_{0}^{t}\int_{0}^{\infty}z\tilde{N}(ds,dz),\eeqnn
where $N(ds,dz)$ is a poisson random measure with intensity
$ds\mu(dz)$. By Theorem \ref{t6}, there is a unique strong solution
to \beqnn
x(t)=x(0)+\int_{0}^{t}\sigma(x(s))dB(s)+\int_{0}^{t}\int_{0}^{\infty}\phi(x(s-))z\tilde{N}(ds,dz)+\int_{0}^{t}(b(s)+\beta
x(s))ds,\eeqnn which is just another form of (\ref{4.1}). The
conclusion holds immediately. \qed
\bigskip\bigskip

\noindent\textbf{Acknowledgements}\\

I would like to give my sincere thanks to my supervisor Professor
Zenghu Li for his encouragement and helpful suggestions. This
research is supported by National Science Foundation of China (No.
10721091).

\end{document}